\documentclass[a4paper,11pt, reqno]{amsart}

\usepackage[margin=2.75cm]{geometry}

\usepackage{hyperref}
\usepackage{amsmath,amsfonts}
\usepackage{amsthm}
\usepackage{amssymb}
\usepackage{url}
\usepackage{xcolor}
\usepackage[normalem]{ulem}
\usepackage{adjustbox}
\usepackage{multirow}
\usepackage{tikz}
\usetikzlibrary{shapes,snakes}
\usepackage{pgfplots}
\usetikzlibrary{positioning}

\def\R{\mathbb{R}}
\def\d{\mathrm{d}}

\def\L{\mathcal{L}}

\def\d{\mathrm{D}}

\def\x{\bold{x}}
\def\RR{\mathcal{R}}
\def\tT{\widetilde{T}}
\def\tV{\widetilde{V}}
\def\tE{\widetilde{A}}

\def\tL{\widetilde{\mathcal{L}}}

\newcommand{\dsum}{\displaystyle\sum}

\newtheorem{theorem}{Theorem}

\newtheorem{lemma}[theorem]{Lemma}

\pgfplotsset{compat=newest}
\usetikzlibrary{arrows.meta}
\usetikzlibrary{decorations.markings}
\usepgfplotslibrary{statistics}

\definecolor{gold}{RGB}{255, 215, 0}
\definecolor{royalblue}{RGB}{65, 105, 225}
\definecolor{limegreen}{RGB}{50, 205, 50}
\definecolor{tomato}{RGB}{255, 99, 71}

\allowdisplaybreaks

\let\origmaketitle\maketitle
\def\maketitle{
	\begingroup
	\def\uppercasenonmath##1{} % this disables uppercasing title
	\let\MakeUppercase\relax % this disables uppercasing authors
	\origmaketitle
	\endgroup
}

\begin{document}

	\title[]{\Large Fixed Topology Minimum-Length Trees with Neighborhoods}

	\author[V. Blanco, G. Gonz\'alez, \MakeLowercase{and} J. Puerto]{{\large V\'ictor Blanco$^\dagger$, Gabriel Gonz\'alez$^\dagger$, and Justo Puerto$^\ddagger$}\medskip\\
		$^\dagger$Institute of Mathematics (IMAG), Universidad de Granada\\
		$^\ddagger$Institute of Mathematics (IMUS), Universidad de Sevilla\\
		\texttt{vblanco@ugr.es}, \texttt{ggdominguez@ugr.es}, \texttt{puerto@us.es}
	}
	
	\date{\today}
	
	\maketitle

\begin{abstract}
In this paper, we introduce the Fixed Topology Minimum-Length Tree with Neighborhood Problem, which aims to embed a rooted tree-shaped graph into a $d$-dimensional metric space while minimizing its total length provided that the nodes must be embedded to some restricted areas. This problem has significant applications in efficiently routing cables or pipelines in engineering designs. We propose novel mathematical optimization-based approaches to solve different versions of the problem based on the domain for the embedding. In cases where the embedding maps to a continuous space, we provide several Mixed Integer Nonlinear Optimization formulations. If the embedding is to a network, we derive a mixed integer linear programming formulation as well as a dimensionality reduction methodology that allows for solving larger problems in less CPU time. A data-driven methodology is also proposed to construct a proper network based on the instance of the problem. We report the results of a battery of computational experiments that validate our proposal.
\end{abstract}

%\MSC[2010] 90C27 \sep  05C05 \sep 90C11 \sep 90C30.
%\end{keyword}

\section*{Introduction \label{sec:introduction}}

Determining the routes of cables or pipelines represents one of the most challenging phases when designing engineering structures. These designs arise when constructing  ships, buildings, or bridges, where designers must strategically trace the different cables throughout intricate spaces, avoiding incompatible interactions between them, maintaining safety distances, and meeting other complex technical requirements. In these types of infrastructures, it is common to be provided with a reference diagram that serves as a visual layout of certain connections between various elements linking the cables. This \emph{schematic planning} explicitly indicates how the different elements that are part of the system are linked and how they should be placed. These connections must be followed to ensure the system's proper performance and represent the  topological structure that the cables must adhere to in order to connect the different elements requiring cable service, such as valves and outlets. While this representation describes the connections between the different elements, it does not specify the exact positions of the elements in the embedded space. Instead, each of these elements has to be placed in certain given areas. This flexibility in tracing the cables in the workspace implies deciding among infinitely many routes for them.

In this paper, we analyze the combinatorial optimization problem behind this type of designs. Specifically, we introduce the so-called \emph{Fixed Topology Minimum-Length Tree with Neighborhood Problem (FTSTNP)}. In this problem, we are given a set of regions where the different elements are to be placed, as well as a tree-shaped hierarchical structure of the different relationships between these elements that inform about all the links between them. The goal is to embed this tree into the designated space, placing the elements within the specified regions in a manner that minimizes the total length of the tracing.

The problem shares some characteristics with various Combinatorial Optimization problems with Neighborhoods analyzed in the literature (see e.g., \cite{blanco2017_MSTN,blanco2022_HLPN,puerto2024hampered,disser2014rectilinear,marcucci2024shortest,espejo2022minimum},
among others). In these types of problems, the goal is to embed a combinatorial structure (such as a path, a spanning tree, a hub network, or a matching) into a given $d$-dimensional space. This embedding requires that the embedded nodes belong to predefined regions (also known as \emph{neighborhoods}) and that a specified metric of the embedded graph is optimized. These neighborhoods represent either some kind of uncertainty in the actual positions of the nodes or the flexibility that users allow for the positions of the elements that want to be located in the embedding space.

In particular, in the Minimum Spanning Tree with Neighborhood Problem (MSTNP), a set of neighborhoods is provided and the goal is to construct a minimum-distance spanning tree where the node positions are determined within the neighborhoods~\cite{blanco2017_MSTN,yang2007minimum,dorrigiv2015minimum,disser2014rectilinear}. It is already known that neighborhood versions of  combinatorial optimization problems which are known to be solvable in polynomial-time, become NP-hard when combined with the requirement of continuous location within different neighborhoods.

Although the FTSTNP analyzed in this paper can be seen as a neighborhood version of a combinatorial optimization problem ---specifically, embedding a tree in a given space--- the goal is not to minimize the sum of each edge's length, but rather the overall length of the cables routed through this structure allowing for extra connections points provided that they reduce the overall length. This problem is closely related to a neighborhood version of the Steiner Tree Problem (STP).

The (continuous) STP involves finding a minimum length (for a given metric) spanning tree on a graph that includes at least a specified set of nodes (known  as \emph{terminals}), with additional nodes also allowed (known as \emph{Steiner points}). The STP is an NP-hard problem, even in the Euclidean case~\cite{fampa2016overview,garey1977complexity}. It has garnered significant attention due to both its theoretical challenges and practical implications.  Specifically, in continuous spaces, when given a set of points, the STP can be seen as the problem of constructing the shortest network connecting the points in the set according to some distance measure. Thus, the STP serves as a generalization of the Weber problem, where the goal is to identify the point that minimizes the sum of distances to a given set of points~\cite{drezner2002weber,fekete2005continuous}. Furthermore, the STP can be seen as a generalization of the Weber problem~\cite{weber1929standort}, which is one of the fundamental problems in Facility Location.

In this paper, we propose a mathematical optimization-based approach for the FTSTNP, which we derive by showing its equivalence to a particular neighborhood version of the STP. Thus, the FTSTNP also shares some characteristics with the so-called Fixed Topology Steiner Tree (FTST) problem~\cite{sankoff1975locating,kiefner2016minimizing,lusheng1999fixed}. This problem consists of embedding the Steiner nodes of a given Steiner Tree such that the overall length of the tree (with the given topology) is minimized. However, in the FTST problem, the underlying tree is known and fixed, whereas in our approach we only have precedence relations between the nodes, but the input graph derived in our construction is no longer a tree. Furthermore, as far as we know, the FTST problem has only been analyzed with rectilinear distance used to measure lengths, whereas in our case, we provide an approach for general $\ell_\tau$ distances in $d$-dimensional spaces. Additionally, the previous approaches are based on combinatorial constructions and analyze the theoretical complexity of the problems. We propose an exact algorithmic framework based on mixed integer linear and nonlinear optimization formulations that can be implemented using available off-the-shelf optimization software, serving as a decision aid tool for practitioners.

As already mentioned, the FTSTNP provides a solution to problems arising in the design of cables-and-pipelines routing. Several papers have addressed this challenging problem, which finds applications in naval, wind farms, and circuit design contexts (see e.g. ~\cite{blanco2021network,blanco2023pipelines,fischetti2018optimizing}). Thus, our approach has a direct impact on challenging industrial problems that arise in many different fields. With this in mind, we provide two distinct frameworks for the FTSTNP based on the domain of the solutions. On the one hand, we analyze the continuous-domain problem, where the embedding of the tree is allowed in the entire $d$-dimensional space $\R^d$. The nodes of the embedded tree will be allowed to be located within the given neighborhoods. The shapes of the neighborhoods will be general and assumed to be unions of second-order cone representable sets. In the second framework that we study, the discrete-domain problem, an embedded graph in $\R^d$ is assumed to be given, and the tree is assumed to be embedded within this graph. The embedding should map the nodes of the tree to the nodes of the graph and the edges to shortest paths in the graph linking the extreme embedded nodes. Although the latter approach may seem more restrictive (not all the positions in $\R^d$ are allowed in the embeddings), it allows for the incorporation of different technical requirements (such as obstacles or less preferred links) that would be difficult to consider in the continuous-domain approach.

Summarizing, the main contributions of this paper are:
\begin{itemize}
\item We propose a novel mathematical optimization-based approach to solve the FTSTNP under different domains that are useful for practical applications: continuous and discrete.
\item For the continuous case, we provide a general Mixed Integer Second Order Cone model. Additionally, for the Euclidean and rectilinear cases, we offer an alternative, more specialized formulation for the problem. 
\item For the discrete case, we derive a Mixed Integer Linear Programming model that allows further generalizations, enabling designers to create more accurate models customized to the specific characteristics of their applications. Additionally, we develop a method to reduce the size of the mathematical optimization model, which directly impacts the size of the problem that optimization solvers can handle.
\item We report the results of an extensive battery of experiments to demonstrate the validity of our approach.
\end{itemize}

The rest of the paper is organized into five sections. In Section \ref{sec:problem} we introduce the FTTNP. Section \ref{sec:conti} is devoted to demonstrating the equivalence of the problem to a particular version of the Fixed Topology Minimum-Length Steiner Tree with Neighborhood Problem, which is constructed from an instance of our problem.  In Section \ref{sec:conti}, we develop a general model for the continuous case, which allows the use of general metrics in any $d$-dimensional space. The Euclidean and rectilinear cases are further analyzed to strengthen the mathematical formulations. On the other hand, Section \ref{sec:discrete} analyzes the discrete case on general graph embeddings and details the mathematical model for this framework.  In Section \ref{sec:experiments}, we report the results of our computational experiments, analyzing the computational performance of our algorithms as well as the differences of the solutions obtained with the different frameworks. Finally, in Section \ref{sec:conclusions}, we draw conclusions and outline future lines of research.
% -------------------------------------------
% -------------------------------------------
% -------------------------------------------

\section{Fixed Topology Minimum-Length Trees with Neighborhoods \label{sec:problem}}

In this section we introduce the problem under analysis and the notation that will be used throughout this paper.

We are given a rooted tree-shaped graph $T=(V,A)$, where $V=\{v_{0}, \ldots, v_{n}\}$ (here, $v_{0}$ is identified as the root of the tree), representing the set of nodes, and $A$ denotes the set of arcs. We denote by $N=\{0, \ldots, n\}$ the index set for the nodes. The hierarchical (parent-child) relations are provided by the mapping $\varrho: V\backslash{v_0} \rightarrow V$ (where $\varrho(v)$ is defined as the unique parent of node $v\in V\backslash\{v_0\}$). To simplify, we use $\varrho^{-1}$ to represent the inverse set mapping, which maps each non-leaf node to its set of children. In other words, $\varrho^{-1}: V \rightarrow \L^{2^V}$, such that $\varrho^{-1}(v) = \{w \in V: \varrho(w) = v\}$. Note that this topology only specifies the connections between the different nodes, but not the positioning of the nodes and therefore nor the lengths of the arcs.

In addition, each node in $V$ is endowed with a neighborhood. For each $i \in N$, we are given a region $\mathcal{R}_{i} \subset \R^{d}$ indicating the positions in $\R^d$ where node $v_i$ can be embedded. Finally, we are given a distance, $\d$ in $\R^d$, to measure lengths in the embedding. 

The goal of the \emph{Fixed Topology Minimum-Length Trees with Neighborhoods Problem} (FTTNP, for short) is to find an embedding $\rho: V \rightarrow \R^d$ such that $\rho(v_i) \in \mathcal{R}_i$, it fulfills the hierarchies of the tree $T$ and minimizes the overall length of the embedded tree.

\begin{figure}[h] 
\begin{center}
\includegraphics[trim={1.7cm 0.5cm 0.4cm 0.4cm},clip,width=1\textwidth]{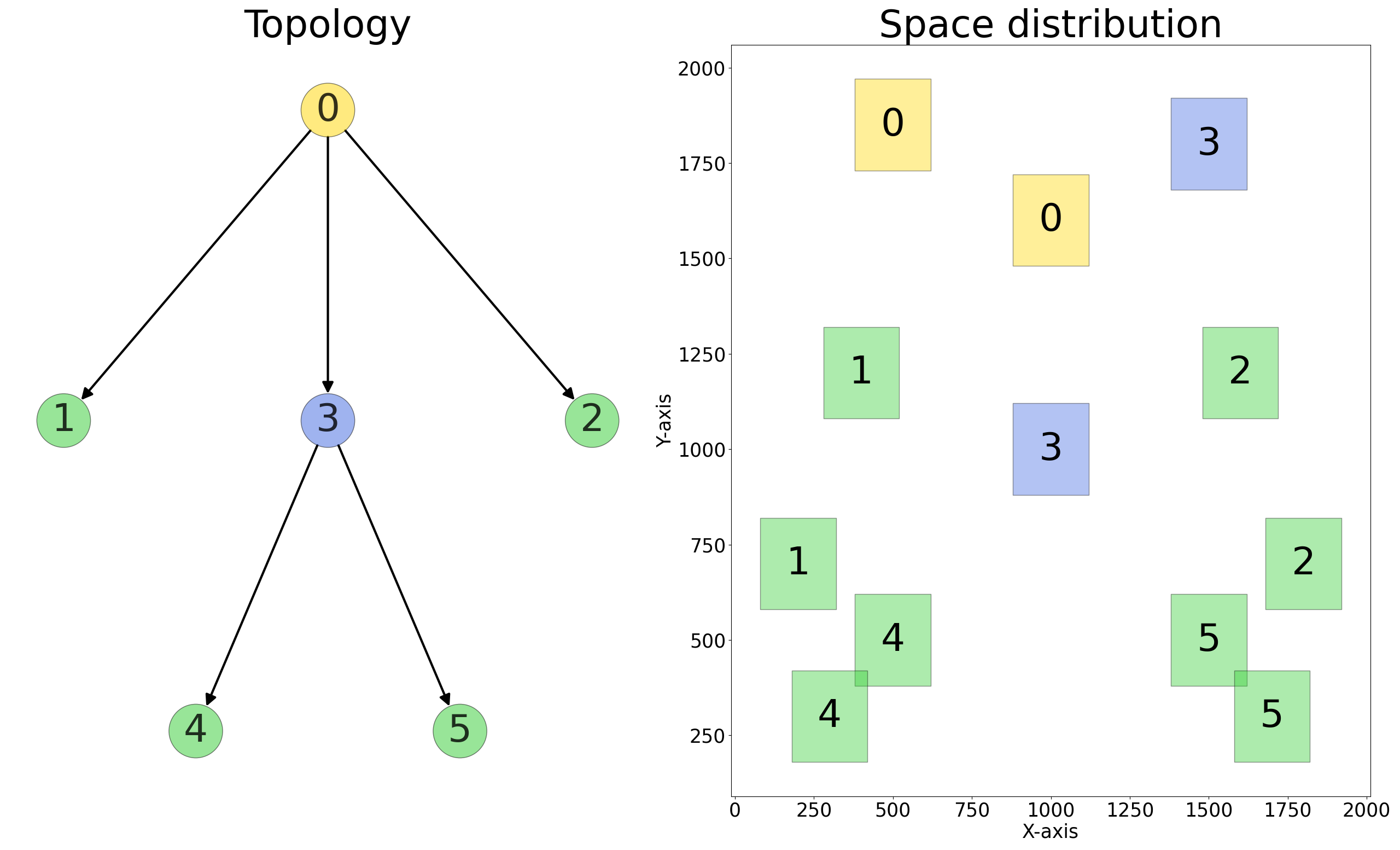}
\end{center}
\caption{Representation of the graph $T$ and the spaces assigned to each node.\label{fig:totally_unsolved}}
\end{figure}

In Figure \ref{fig:totally_unsolved} we show an illustrative instance for the FTTNP. In the left plot we show the rooted tree $T$ indicating the topology that must be followed. Note that by following these requirements, in the embedding to be found, we will account exclusively for the lengths of the paths linking the children to their parents in that tree. In the right picture we show the neighborhoods for each of the nodes of the tree. In the embedding, the root node can only be positioned in the space determined by the $0$-labeled rectangles (that is, in one of these rectangles), node $1$ in the $1$-labeled regions, etc. 

As mentioned in the introduction of this paper, one of the the most practical applications of this problem comes from cable routing. In the design of a building or a ship, once the main pipelines or cable trunks have been laid out in a 3D space according to constructibility criteria (for further details, see \cite{blanco2021network,blanco2023pipelines}), the challenge arises of optimally routing cables from the primary pipelines to terminals while minimizing cable length. Guided by a schematic plan, valves are strategically positioned within designated areas to fulfill technical requirements. This challenge is contextualized within FTTNP, where the schematic of each pipeline forms a rooted tree structure, with valves optimizing cable routes between pipelines, terminals, and other valves.

Finally, we would like to highlight that the FTTNP is not a Fixed Topology Minimum-Length Spanning Tree Problem with Neighborhoods (MSTN). 

Note that the (fixed topology) MSTN is mathematically formulated as:
\begin{align}
\min &\ \sum_{(v_j, v_k) \in A}  \d(\x_j,\x_k) \tag{FT-MSTN}\label{fttnp}\\
\mbox{s.t. } & \x_j \in \mathcal{R}_j, \forall j\in N.\nonumber
\end{align}
That is, each node of the given tree ($v_j$) is embedded within its neighborhood (at coordinates $\x_j\in \mathcal{R}_j$) minimizing the overall sum of the arc lengths of the tree. Nevertheless no extra intermediate points can be added to minimize the overall length.

The problem described above can be efficiently solved (in polynomial time) by most off-the-shelf software  if the distance measure is an $\ell_\tau$-norm or a polyhedral norm. This is because in those cases, the problem can be reformulated as a second-order cone or linear optimization problem (see e.g. \cite{blanco2014revisiting,blanco2024minimal}) and solved using interior point methods \cite{kuo2004interior}. However, this approach may not yield desirable solutions for practical applications. Specifically, since the distances considered in the objective function of \eqref{fttnp} are calculated edge-by-edge, there is no encouragement for cases where multiple edges share part of the route from the root node to the leaves, which could result in less costly solutions. Moreover, solving \eqref{fttnp} may lead to solutions where different geometric edges overlap or cross, resulting in unrealistic outcomes for practical applications, such as those encountered in engineering designs where the edges may represent cables. The aim of this paper is to develop an alternative methodology to obtain realistic solutions for the problem, taking into account the various scenarios that should be promoted or avoided when constructing solutions to be implemented in real-world situations.

\begin{figure}[h]
\begin{center}
\begin{tikzpicture}[scale=0.9]
    \node (A1) at (0,0) [circle, draw, fill=gold, inner sep=1, opacity=.5, text opacity=1, minimum size=6mm] {\normalsize $0$};
    \node (A2) at (-1,-3) [circle, draw, fill=limegreen, inner sep=1, opacity=.5, text opacity=1, minimum size=6mm] {\normalsize $1$};
    \node (A3) at (1,-3) [circle, draw, fill=limegreen, inner sep=1, opacity=.5, text opacity=1, minimum size=6mm] {\normalsize $2$};
    \draw[-Latex, line width=0.3mm] (A1)--(A2);
    \draw[-Latex, line width=0.3mm] (A1)--(A3);
\end{tikzpicture}~\hspace*{0.3cm}
% ------------------------------------
\begin{tikzpicture}[scale=0.9]
\draw[draw=gold, fill=gold, opacity=0.5] (0,0) rectangle ++(1,1);
\draw[draw=limegreen, fill=limegreen, opacity=0.5] (-1,-3) rectangle ++(1,1);
\draw[draw=limegreen, fill=limegreen, opacity=0.5] (1,-3) rectangle ++(1,1);
\node (A0) at (0.5,0.5) {$0$};
\node (A1) at (-0.5,-2.5) {$1$};
\node (A2) at (1.5,-2.5) {$2$};
\end{tikzpicture}~\hspace*{0.3cm}
% ------------------------------------
\begin{tikzpicture}[scale=0.9]
\draw[draw=gold, fill=gold, opacity=0.5] (0,0) rectangle ++(1,1);
\draw[draw=limegreen, fill=limegreen, opacity=0.5] (-1,-3) rectangle ++(1,1);
\draw[draw=limegreen, fill=limegreen, opacity=0.5] (1,-3) rectangle ++(1,1);
    \node (A1) at (0.5,0) [inner sep=0, color=tomato]  {\small $\bullet$};
    \node (A2) at (0,-2) [inner sep=0, color=tomato]{\normalsize $\bullet$};
    \node (A3) at (1,-2) [inner sep=0, color=tomato] {\normalsize $\bullet$};
\draw[-Latex, line width=0.3mm] (A1)--(A2);
\draw[-Latex, line width=0.3mm] (A1)--(A3);
\node (A0) at (0.5,0.5) {$0$};
\node (A1) at (-0.5,-2.5) {$1$};
\node (A2) at (1.5,-2.5) {$2$};
\end{tikzpicture}~\hspace*{0.3cm}
% ------------------------------------
\begin{tikzpicture}[scale=0.9]
\draw[draw=gold, fill=gold, opacity=0.5] (0,0) rectangle ++(1,1);
\draw[draw=limegreen, fill=limegreen, opacity=0.5] (-1,-3) rectangle ++(1,1);
\draw[draw=limegreen, fill=limegreen, opacity=0.5] (1,-3) rectangle ++(1,1);
    \node (A1) at (0.5,0) [inner sep=0, color=tomato]  {\normalsize $\bullet$};
    \node (A2) at (0,-2) [inner sep=0, color=tomato]{\normalsize $\bullet$};
    \node (A3) at (1,-2) [inner sep=0, color=tomato] {\normalsize $\bullet$};
\draw[-, line width=0.3mm] (A1)--(0.5,-2);
\draw[-Latex, line width=0.3mm] (0.5,-2)--(A2);
\draw[-Latex, line width=0.3mm] (0.5,-2)--(A3);
\node (A0) at (0.5,0.5) {$0$};
\node (A1) at (-0.5,-2.5) {$1$};
\node (A2) at (1.5,-2.5) {$2$};
\end{tikzpicture}
\end{center}
\caption{A toy instance to show the distance between the MSTN and FTTNP.\label{mstn}}
\end{figure}
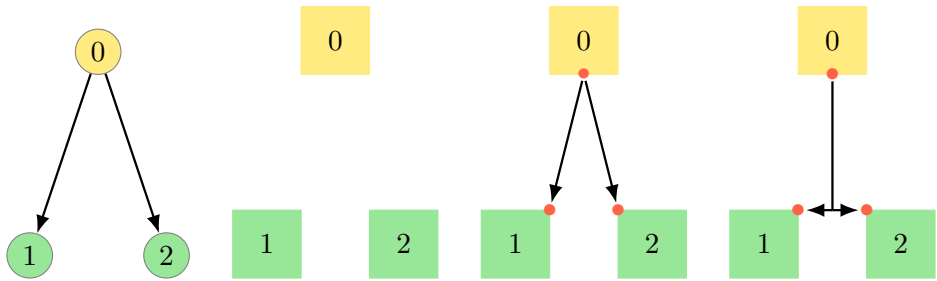

In Figure \ref{mstn} we illustrate the solution of this problem compared to the FTTNP. The two left plots represent a toy instance for the problem (the tree topology and the neighborhoods, respectively). The third picture shows the solution of the FT-MSTNP. That is, the nodes are embedded in their neighborhoods such that the sum of the lengths of the two embedded arcs is minimized. In contrast, the far-right picture shows a solution of the FTTNP, which, instead of separately minimizing each of the lengths, minimizes the design of the overall trace of the system. 

Here, we present an alternative methodology for embedding a given tree into a $d$-dimensional space while minimizing the space used by the arcs of the tree in the embedding. Specifically, we propose methodologies based on constructing minimum-length Steiner Trees with specifications provided by the input tree. Along these lines, we develop two approaches, each with its own advantages and disadvantages, based on whether or not the positions of the embedding are restricted to a finite set of positions in the neighborhoods. Both approaches differ in the domain of the embedding, and can be useful in different situations. Furthermore, we provide a data-driven methodology to derive approximate solutions for the continuous case by using the discrete approach on a certain graph obtained from the original topology of the tree and neighborhoods.

\section{The continuous-domain FTTNP \label{sec:conti}}

In this section, we extend \eqref{fttnp} to allow for the sharing of paths from the root node to the leaves among certain edges, aiming to reduce the overall length of the tree. We achieve this by formulating the FTTNP. Our methodology has two main components: (1) transforming the input tree, $T$, by adding an additional set of nodes and edges that incorporate the overall design length into the problem; and (2) determining the tree structure and the geometric positions of the nodes in $\R^d$, with the possibility of discarding some of the newly introduced sets of nodes, using a Steiner tree-like model.

In the following subsections, we detail these steps.

\subsection{Graph Transformation \label{subsec:graph}}

The first step of our approach involves constructing a new graph from the original input tree that allows us to adequately model the FTTPN as a suitable mathematical optimization problem. The key idea of this transformation is to enable the paths resulting from embedding the arcs in the input graph to be shared among different arcs. Thus, we replace parent-children arcs by two-legs paths with auxiliary nodes which are all connected between them. In the following, we detail this construction.

Given the input rooted tree $T=(V,A)$, we construct a new graph $\tT=(\tV,\tE)$ as follows:

% For each $v \in V$with $|\varrho^{-1}(v)| \geq 2$ (indicating that node $v$ has at least two children in the tree),
\begin{itemize}
    \item We initialize $\tT=T$.%, with $\tV=V$ and $\tE=A$.
    \item For each $v\in V$ with $|\varrho^{-1}(v)| \geq 2$ (indicating that node $v$ has at least two children in the tree), we remove all arcs between $v$ and its children, i.e.,
$$
\tE \leftarrow \tE \Big\backslash\Big\{(v,w): w \in p^{-1}(v)\Big\}, \text{for all $v\in V$.}
$$
    \item For each node $v\in \tV$, $|\varrho^{-1}(v)|-1$  new nodes, $w_1^v, \ldots, w_{n_v-1}^v$, are added to $\tV$, as well as arcs with source $v$ and target to each of these new nodes, and arcs from each of these new nodes to each of the children:% linking thoses nodes with the nodes in $\tV$. Specifically, for each auxiliary node, $1+|p^{-1}(v)|$ new arcs are introduced as follows:
    $$
    \tE \leftarrow \tE \cup \Big\{(v,w_j^v)\Big\} \cup \Big\{(w_j^v,w): w \in \varrho^{-1}(v)\Big\},
    $$
    for all $j=1, \ldots, n_{v}-1$.
    Furthermore, all the auxiliary nodes for node $v$ are completely linked in the new graph, i.e.,
    $$
    \tE \leftarrow \tE \cup \Big\{(w^v_j,w^v_k): j, k =1, \ldots, |\varrho^{-1}(v)|-1, j\neq k\Big\}.
    $$
\end{itemize}
Note that with this construction the graph $\tT$ is no longer a tree, but it is connected.

\begin{figure}[h] 
\begin{center}
\includegraphics[trim={1.7cm 0.5cm 0.4cm 0.4cm},clip,width=1\textwidth]{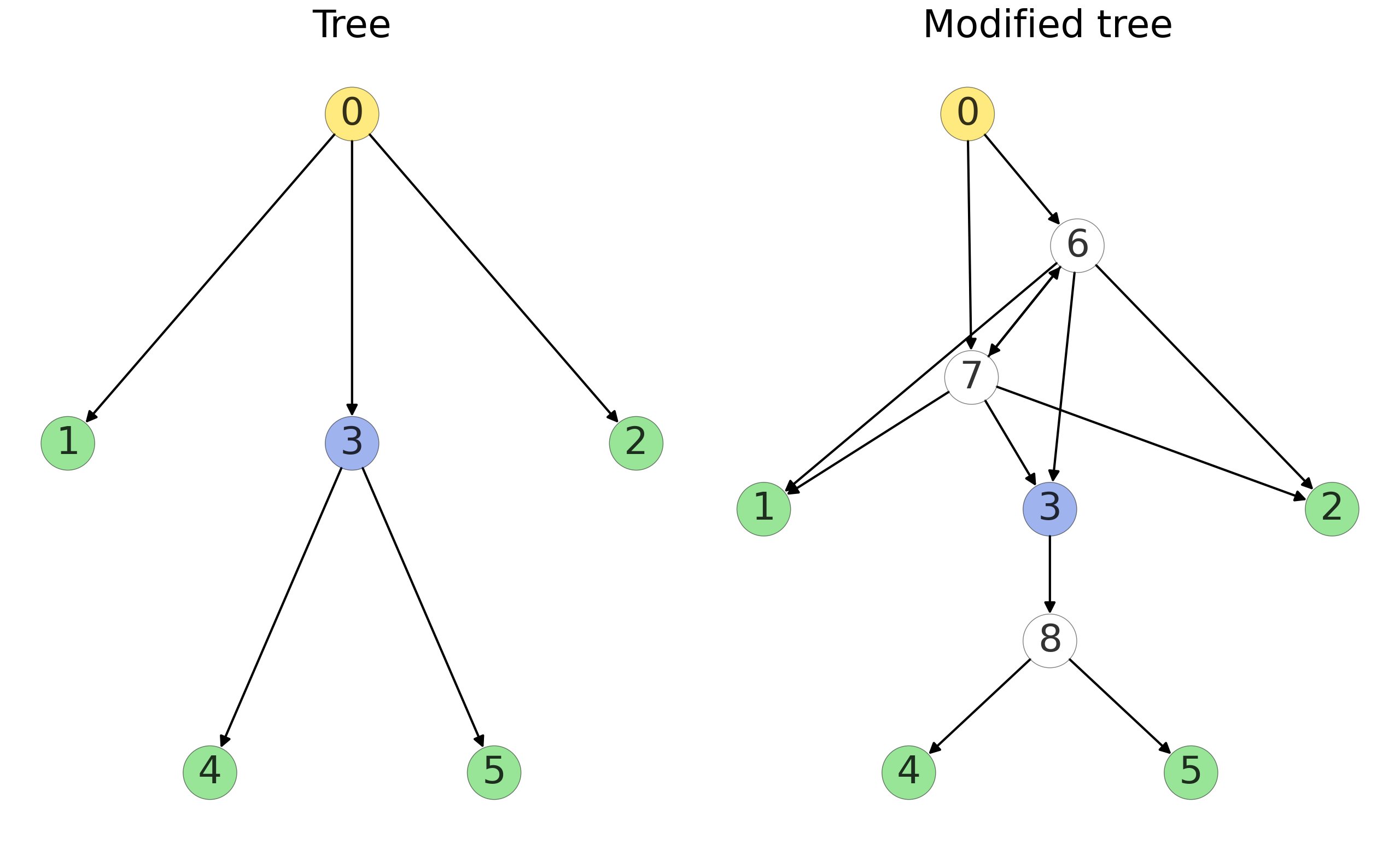}
\end{center}
\caption{Representation of tree graph $T$ and the transformation $\tT$.\label{fig:ejem_trans}}
\end{figure}

In Figure \ref{fig:ejem_trans}, we illustrate this transformation. On the left side, a simple tree topology $T$ is shown. The right picture depicts the transformed graph $\tT$ for this tree. Note that the spatial positions of these nodes are not representative of the embedding we aim to construct in our approach; only the topology of the tree. The new graph includes all the nodes of $T$, as well as the new nodes $6$, $7$, and $8$. Nodes $6$ and $7$ are incorporated because node $0$ has three children (nodes $1$, $2$, and $3$). Thus, arcs $(0,6)$ and $(0, 7)$ are added, as well as arcs linking these nodes with the children nodes. Similarly, node $8$ appears because node $3$ has two children, but only the arcs linking the parent to it, and it to the children, are added.

Furthermore, each node $v\in V$ in the tree $T$ is endowed with a neighborhood $\RR_v \subseteq \R^d$. In the transformed graph $\tT =(\tV,\tE)$, the nodes in $V$ retain their original neighborhoods. However, for nodes $w \in \tV\backslash V$, the neighborhoods are defined as $\RR_w = \R^d$, indicating that any position in space is feasible for these nodes.

\subsection{The Steiner tree problem on $\tT$}

In this phase, we equivalently model the FTTNP as a restricted Steiner tree problem over the transformed graph $\tT$. 

Given a (di)graph and a given set of terminal nodes from the graph's node set, the goal of the Steiner tree problem is to construct a minimum length tree in the graph that at least uses the terminal nodes, while also being allowed to use some of the other nodes from the graph. 

A Minimum Length Spanning Trees connects the nodes for which there is a direct arc linking them directly, which might not always be the most cost-effective way to lay cables if intermediate nodes (additional points where cables can be laid) can help to reduce the total cable length. In contrast Minimum Length Steiner Trees allow for the use of additional nodes that can significantly reduce the total length of the resulting network. Furthermore, Steiner trees provide more flexibility in the placement of cables, as they are not constrained to a fixed number of nodes. This can be particularly useful in physical spaces where laying cable directly between all original points is not feasible or too costly. In any case, if no new nodes are needed, the Steiner tree problem will decide not to incorporate them, and its solution will coincide with the minimum spanning tree.

Denoting by $\tV = \{v_0, \ldots, v_{n+1},v_{n+2}, \ldots, v_{|\tV|}\}$ the labels for the nodes of the transformed graph, and by $N' = N \cup \{n+1, \ldots, |\tV|\}$ its index set, the Steiner tree problem that we analyze here can be formally formulated as follows:

\begin{align}
    \min_{S \subset \tV\atop B \subset \tE} &\dsum_{(v_j,v_k) \in B} \d(\x_j,\x_k)\\
    \mbox{s.t. } & v_j \in S, \label{st:1}\\
    & (S,B) \mbox{ is a tree},\label{st:2}\\
    & \x_j \in \RR_j, j \in N,\label{st:3}\\
    & \x_j \in \R^d, j \in N'\backslash N,\label{st:4}
\end{align}
where the coordinates $\x_j$ indicates the embedding of the nodes in $\R^d$.

The objective function minimizes the total length of the tree. Constraint \ref{st:1} ensures that the \emph{terminals} (original nodes of the tree) are in the resulting subgraph (which is required to be a tree by constraint \eqref{st:2}). Constraints \eqref{st:3} and \eqref{st:4} define the domains for the embedding. Here, abusing of notation we identify $\RR_{v_j}$ with $\RR_j$. For terminal nodes, the embedding is constrained to their neighborhoods, and for the Steiner points, it is allowed in the whole space.

Observe that the problem above is not exactly a Steiner tree problem on a Euclidean space since some hierarchical conditions are imposed to the problem to respect the required topology. 

\subsection{Mathematical Optimization Model \label{subsec:mathmodel}}

To construct a model for the FTSTNP, we formulate it as a Mixed Integer Nonlinear Optimization problem through a flow-based formulation. Apart from the $\x$-variables already described, we use the following sets of variables: 
\begin{align*}
    y_{jk} = &\begin{cases}
    1, & \mbox{if arc $(v_j,v_k)$ is part of the resulting tree,}\\
    0, & \mbox{otherwise,}
\end{cases}\\
&\mbox{for all } (v_j,v_k) \in \tE.\\
% ----------------------
f^\ell_{jk} = &\begin{cases}
    1, & \mbox{if arc $(v_j,v_k)$ is used in the path from the root node to the leaf node $v_\ell$,}\\
    0, & \mbox{otherwise,}
\end{cases}\\
&\mbox{for all $(v_j,v_k) \in \tE$, and $v_\ell \in \tL$ (here $\tL$ denotes the set of leaves in the original tree).}\\
% ----------------------
d_{jk}\in &\R_{+} \mbox{: Length of arc linking $\x_j$ and $\x_{k}$, for all $(v_j, v_k) \in \tE$.}
\end{align*}

The model can be stated as follows:
\begin{align}
\min & \sum_{(j, k) \in \tE} d_{jk} y_{jk} \label{m1:obj}\\
\text{s.t.} & \sum_{k \in N:\atop (v_0,v_k) \in \tE} f_{0k}^{\ell} = 1, &\forall \ell \in \L,\label{m1:ctr3}\\
& \sum_{k \in N:\atop (v_j,v_k) \in \tE} f_{jk}^{\ell} - \sum_{k \in N:\atop (v_k,v_j) \in \tE} f_{kj}^{\ell}=0, &\forall \ell \in \L, j \in N'\backslash N,\label{m1:ctr4}\\
& \sum_{j \in N:\atop (v_j,v_\ell) \in \tE} f_{j \ell}^{\ell} = 1, &\forall \ell \in \L,\label{m1:ctr5}\\
& f_{jk}^{\ell} \leq y_{jk}, &\forall \ell \in \L, j, k \in N', (v_j,v_k) \in \tE,\label{m1:ctr6}\\
& d_{jk} \geq \d(\x_j,\x_{k}), &\forall j, k \in N': (v_j, v_k) \in \tE,\label{m1:ctr7}\\
& \x_j\in \RR_j, &\forall j \in N, \label{m1:ctr8}\\
& \x_j\in \R^d, &\forall j \in N'\backslash N, \label{m1:ctr9}\\
& f_{jk}^\ell \geq 0, &\forall j, k \in N': (v_j, v_k) \in \tE, \ell \in \L,\\
& y_{jk} \in \{0,1\}, &\forall j, k \in N'.
\end{align}
In the above formulation, the objective function \eqref{m1:obj} aims to minimize the total length of the activated arcs of the tree. Constraints \eqref{m1:ctr3} ensure that exactly an arc is activated in the solution, departing from the root. Constraints \eqref{m1:ctr4} represent the flow conservation constraint for each leaf node. Constraints \eqref{m1:ctr5} guarantee that each leaf node is reached by exactly a single arc. Constraints \eqref{m1:ctr6} prevent the activation of flows on arcs that are not utilized. Constraints \eqref{m1:ctr7} ensure the proper definition of the lengths of the arcs in the tree. Finally, Constraints \eqref{m1:ctr8} and \eqref{m1:ctr9} are the domains for the embedding (the original nodes in the neighborhoods and the auxiliary nodes freely in $\R^d$). Note that the domain of the $f$-variables can be relaxed to the nonnegative orthant.

Note that the above formulation is a Mixed Integer Nonlinear Optimization problem due to the objective function and constraints \eqref{m1:ctr6} and \eqref{m1:ctr7}. However, the objective function can be linearized exactly by applying McCormick envelopes, and $\ell_p$ or polyhedral norm-based distances are efficiently represented through Second Order Cone Constraints (see e.g., \cite{blanco2014revisiting,blanco2024minimal}) or linear inequalities. Additionally, if the neighborhoods are second-order cone representable or a union of second-order cone representable sets, one can still keep a Mixed Integer Second Order Cone Optimization model, that can be solved, at least for medium-sized instances, using off-the-shelf software such as Gurobi, CPLEX, or FICO. 

\subsubsection{Special Cases: $\ell_{1}$ and $\ell_{2}$ norm based distances \label{subsec:Improved}}

The Steiner Tree problem with particular distance choices for $\d$ has been studied in the literature in order to exploit the geometric properties to be incorporated into the mathematical optimization models (see e.g. \cite{brazil2009novel,Brazil2014, fampa2016overview, du1993minimum, garey1977rectilinear, Hanan1966, Maculan2000, Ouzia2022, Pinto2023}). Although our problem is not strictly a Steiner tree problem on a continuous space (because of the requirements imposed by the topology of the tree), some of the properties of Steiner trees still apply to our problem. 

Specifically, one can project out the flow variables in our model by  replacing Constraints \eqref{m1:ctr3}, \eqref{m1:ctr4}, \eqref{m1:ctr5}, and \eqref{m1:ctr6} in Subsection \ref{subsec:mathmodel} by the following constraints that only involve the $y$-variables:
\begin{align}
& \sum_{k \in N:\atop (v_j,v_k) \in \tE}y_{jk} = 1, & \forall j \in N'\backslash \L, \label{eq:l8}\\
& y_{jk} + y_{kj} \leq 1, & \forall j, k \in N': (v_j,v_k), (v_k,v_j) \in \tE, \label{eq:l10}\\
& \sum_{k\in N':\atop (v_k,v_j) \in \tE} y_{kj} = 1, & \forall j \in N\backslash \{0\}.\label{eq:l11}
\end{align}
The first set of constraints specify that from an original node of $T$, one can only proceed to a new \emph{Steiner} node; the second set of constraints ensures that one is allowed to use an arc in a single direction; the third set of constraints require that only one arc arrives at a node of the original tree, $T$.

Specifically, for $\ell_p$ norms, one can derive conditions on the degrees of the auxiliary nodes in the Steiner tree, as stated in the following result.
\begin{lemma}[\cite{du1993minimum}]
    Let $\d$ be a distance induced by a norm $\|\cdot\|$, i.e., $\d(a,b) = \|a-b\|$, for $a, b \in \R^{d}$. If $\|\cdot\|$ is differentiable and strictly convex, then the degree of $v_j$ for $j\in N'\backslash N$ is exactly $3$.
\end{lemma}
Note that the the above result applies to any $\ell_p$ norm with $1<p<\infty$, particularly to the Euclidean norm. For more information on the approximability of the Steiner tree problem in $\ell_p$-metrics, see \cite{lpSteiner}. For the $\ell_1$-norm, we have the following result.
\begin{lemma}[\cite{Hanan1966}]
    Let $\d$ be a distance induced by the $\ell_1$ norm, i.e., $\d(a,b) = \sum_{k=1}^d |a_k-b_k|$, for $a, b \in \R^{d}$. Then, the degree of $v_j$ for $j\in N'\backslash N$ is either $3$ or $4$.
\end{lemma}
The result above is based on the construction of the Hanan grid, which is obtained by drawing the hyperplanes parallel to the axes crossing each of the terminal nodes of a classical Steiner tree. 

With the above lemmas, we incorporate these conditions into the model as linear inequalities involving the $y$-variables:

\begin{itemize}
\item If $\d(a,b) = \|a-b\|_p$, for $1<p<\infty$:
\begin{align}
& \sum_{j\in \widetilde{V}\atop\left(i, j\right) \in \widetilde{A}} y_{ij} + \sum_{j\in \widetilde{V}\atop\left(j, i\right) \in \widetilde{A}} y_{ji} = 3, && i \in P. \label{eq:con11}
\end{align}
\item If $\d(a,b) = \|a-b\|_1$:
\begin{align}
    & 3 \leq \sum_{j\in \widetilde{V}\atop\left(i, j\right) \in \widetilde{A}} y_{ij} + \sum_{j\in \widetilde{V}\atop\left(j, i\right) \in \widetilde{A}} y_{ji} \leq 4, && i \in P. \label{eq:l9}
\end{align}
\end{itemize}
Observe that this representation not only reduces the number of constraints in the model, but also drastically reduces the number of variables, since the flow variables $f$ are no longer needed.

\begin{figure}[h] 
\begin{center}
\includegraphics[trim={0.2cm 0.2cm 0.2cm 0.2cm},clip,width=1\textwidth]{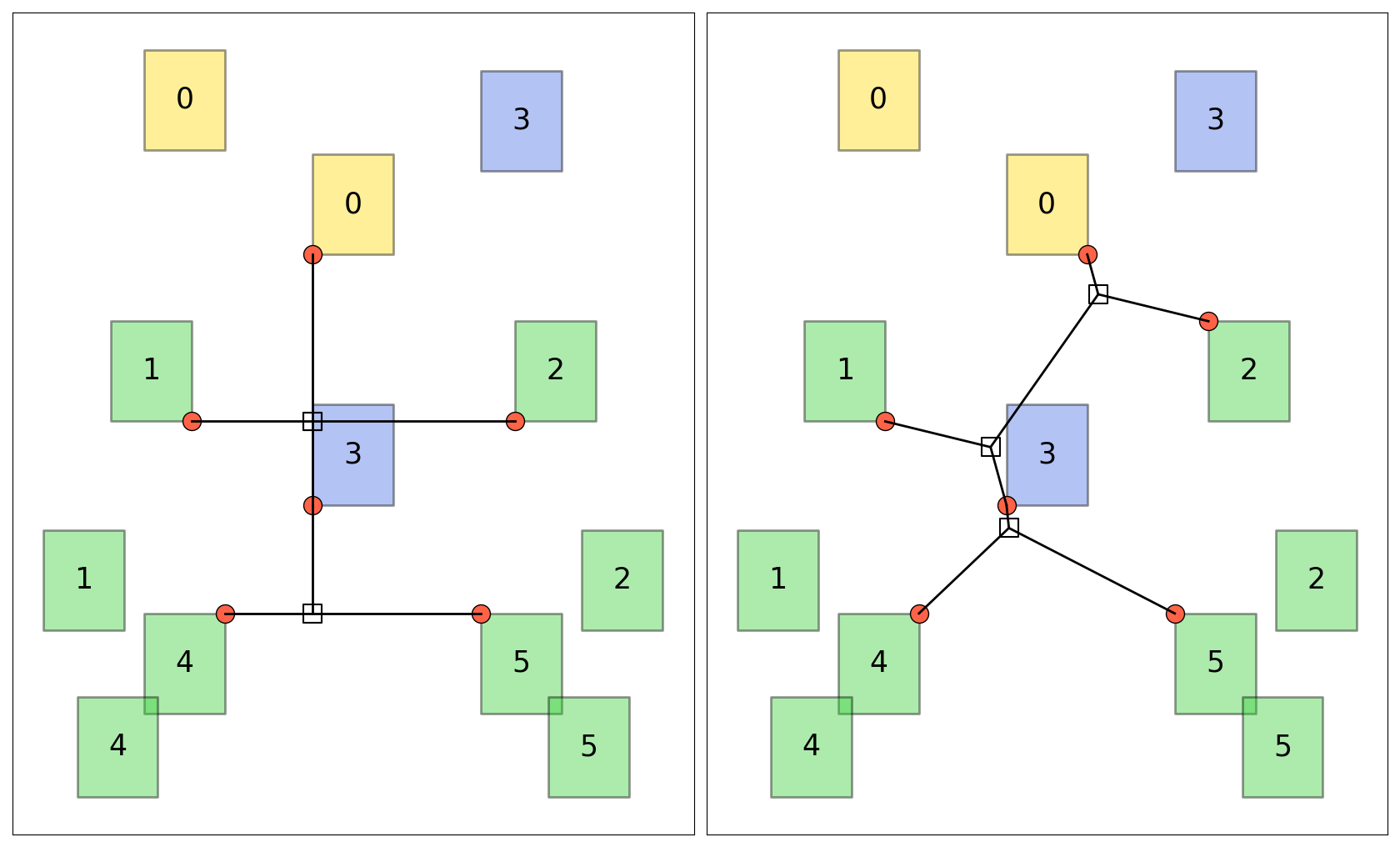}
\end{center}
\caption{Graphical representation of the scenario solution. $\ref{fig:totally_unsolved}$ for the norm $\ell_1$ on the left and $\ell_2$ on the right. \label{fig:xxx}}
\end{figure}

In Figure \ref{fig:xxx} we show the solutions of the same instance that we use to illustrate the problem above, for the $\ell_1$ norm (left) and the $\ell_2$ norm (right).

% --------------------------------------------------
% --------------------------------------------------
% --------------------------------------------------
% --------------------------------------------------
% --------------------------------------------------
% --------------------------------------------------
% --------------------------------------------------
% --------------------------------------------------
% --------------------------------------------------
% --------------------------------------------------
% --------------------------------------------------
% --------------------------------------------------
% --------------------------------------------------
% --------------------------------------------------
% --------------------------------------------------

\section{The discrete-domain FTTNP \label{sec:discrete}}

In the discrete domain case, we assume that we are given an undirected connected graph $G=(U,E)$, with $U=\{u_1, \ldots, u_{|U|}\} \subset \R^d$, and that the embedding maps any node of the tree $T$, $v$, into a node $\rho(v) \in U$. Furthermore, the length of the embedded arc $(\rho(v),\rho(w))$ is the shortest path in $G$ (with any of the available directions in the graph), with source $\rho(v)$ and target $\rho(w)$.

Foundational insights into discrete-based approaches for Steiner problems are provided by \cite{MACULAN1987185, Aneja1980, Diane1993, Aneja2006, Milos2022}, highlighting the relevance of using discrete methods in such optimization tasks.

One of the most popular graphs $G$ that can be used is a grid in $\R^d$ with a certain density and width. If the density and width are large enough, this problem may help to approximate the solution of the continuous-domain case analyzed in the previous section.

Furthermore, this approach has some advantages compared to the previous approach from a practical perspective. On the one hand, the problem can be cast as a mixed integer linear optimization problem that can be more efficiently solved than the nonlinear models derived for the continuous case. On the other hand, the graph approach allows for considering several characteristics in the cable design problem that are not possible (or at least not straightforward) with the proposed approach for the continuous-domain case: 
\begin{enumerate}
    \item Consider obstacles when routing the cables in a given space. This can be easily done by removing from the graph $G$ the nodes and edges that are not allowed to be crossed.
    \item Incorporate preference zones in the routing problem. This can be done by assigning weights to the edges in the graph and taking these weights into account when calculating the shortest path lengths.
    \item Account for (and penalize) the number of direction changes in the routes of the cables. In \cite{blanco2021network,blanco2023pipelines} the authors propose a novel approach to model this technical requirement.
\end{enumerate}

Another advantage of the discrete approach is that it can be useful to approximate the solutions of the continuous-domain versions of the problem. As we will see in our computational experiments, the continuous-domain problems can be very time demanding. The (classical) rectilinear Steiner tree problem in $\R^d$ is known to have a finite dominating set for the Steiner points that can be explicitly constructed using the Hanan grid~\cite{Hanan1966}. Although a finite dominating set is not available for our version of the problem, a similar construction can be done to get good-quality approximations in more reasonable CPU time. We will empirically check the quality of this approximation in Section \ref{sec:experiments}.

In the subsequent subsections, we elaborate on this procedure in detail.

\subsection{Mathematical Optimization Model \label{subsec:math_model}}

The main difference, in terms of the mathematical optimization model, between the continuous and discrete-domain frameworks, is that nonlinearities can be avoided by computing, in preprocessing, the lengths of the arcs involved in the construction, since the nodes are embedded into a finite set of potential points in $\R^d$ (this also helps in defining more sophisticated distances). In contrast, the input data for the discrete problem can be very large (the number of nodes and edges in the graph $G$) which implies a large number of variables and constraints in the model.

Let us start introducing the mathematical model that we propose by describing the variables that will appear in the formulation. For the ease of presentation we introduce the following sets: 
\begin{itemize}
    \item $N_G = \{1, \ldots, |U|\}$: index set for the nodes in $G$. 
    \item $\RR^G_j = \{i \in N_G: u_i \in \RR_j\}$: the nodes of $G$ in the neighborhood of node $j \in N'$ (for those nodes in $N'\backslash N$, $\RR^G_j = N_G$).
\end{itemize}
And the following decision variables:
%$$
%g_{i\ell}^{jk} = \begin{cases}
%    1 & \mbox{if the path from $v_j$ to $v_k$ is routed through edge $\{u_i,u_\ell\}\in E$},\\
%    0 & \mbox{otherwise}
%\end{cases}\; 
%\forall j, k \in N': (v_j,v_k) \in \tV, i, \ell \in N_G: \{u_i,u_\ell\} \in E.
%$$
%$$
%t_{ji} = \begin{cases}
%    1 & \mbox{if  node $v_j\in \tV$ is embeded into node $u_i\in U$},\\
%    0 & \mbox{otherwise}
%\end{cases}\; \forall j \in N', i \in N_G.
%$$
%$$
%q_{i\ell} = \begin{cases}
%    1 & \mbox{if edge $\{u_i,u_\ell\} \in E$ is part of the constructed tree},\\
%    0 & \mbox{otherwise}
%\end{cases}\; \forall  i, \ell \in N_G.
%$$
% ---------
\begin{align*} 
f_{i\ell}^{jk} = &\begin{cases}
    1, & \mbox{if the path from $v_j$ to $v_k$ is routed through edge $\{u_i,u_\ell\}\in E$},\\
    0, & \mbox{otherwise},
\end{cases}\\
&\forall j, k \in N: (v_j,v_k) \in V, i, \ell \in N_G: \{u_i,u_\ell\} \in E.\\
% --
t_{ji} = &\begin{cases}
    1, & \mbox{if  node $v_j\in V$ is embedded into node $u_i\in U$},\\
    0, & \mbox{otherwise},
\end{cases}\\
&\forall j \in N, i \in N_G.\\
% --
q_{i\ell} = &\begin{cases}
    1, & \mbox{if edge $\{u_i,u_\ell\} \in E$ is part of the constructed tree},\\
    0, & \mbox{otherwise},
\end{cases}\\
&\forall  i, \ell \in N_G.
\end{align*}
Note that we do not use a special set of variables to decide which nodes in $V$ are activated as Steiner points or not, but all nodes will be embedded into the new graph. Since the potential Steiner points can be freely embedded in the whole graph, in case a particular node is not used in the final tree, it will coincide with another node, not being accounted for in the objective function and not affecting the length of the tree.

With these variables we propose the following model to solve the discrete-domain FTTNP:

 %\begin{align}
    %\min &\sum_{i, \ell \in N_G:\atop \{u_i,u_\ell\} \in E} d_{i\ell} q_{i\ell} \label{eq:a1}\\
    %\mbox{s. t.} \;
    %&\sum_{i \in \RR_j^G} t_{ji} = 1,  \;\; \forall j \in N', \label{eq:a2}\\
    % &\sum_{i \in N_G:\atop \{u_i,u_\ell\}\in U} \!\!\!\! f^{jk}_{i\ell} -
    % \sum_{i \in N_G:\atop \{u_\ell,u_i\}\in U} \!\!\!\! f^{jk}_{\ell i} = t_{j\ell},
     %\;\;\forall j, k \in N': (v_j, v_k) \in \tE, \ell \in N_G, \label{eq:a3}\\
     %& 
    % \sum_{i \in N_G:\atop \{u_\ell,u_i\}\in U} \!\!\!\!f^{jk}_{\ell i} - \sum_{i \in N_G:\atop \{u_i,u_\ell\}\in U}\!\!\!\! f^{jk}_{i\ell} = t_{k\ell},
    % \;\;\forall j, k \in N': (v_j, v_k) \in \tE, \ell \in N_G, \label{eq:a4}\\
    % & f^{i\ell}_{jk} + f^{\ell i}_{jk} \leq 1,  & \forall i, \ell \in N_G: \{u_i,u_\ell\} \in E, j, k \in N': (v_j,v_k) \in \tE, \label{eq:a6} \\
    % % % -------------------------------------
    %  &q_{i\ell} \geq f^{jk}_{i\ell}, \;\;\forall i, \ell \in N_G: \{u_i,u_\ell\} \in E, j, k \in N': (v_j,v_k) \in \tE, \label{eq:a7} \\
    % % % -------------------------------------
     % &\sum_{\ell \in N_G:\atop \{u_i,u_\ell\}\in E} q_{i\ell} = t_{ji}, \;\; \forall i \in N_G, \label{eq:a8} \\
    % % % -------------------------------------
      %&\sum_{i \in N_G:\atop \{u_i,u_\ell\} \in E} q_{i\ell} \leq 1 - t_{j\ell}, \;\; \forall \ell \in N_G, j\in N'. \label{eq:a9}\\
      %& f_{i\ell}^{jk} \geq 0, \;\;\forall j, k \in N': (v_j, v_k) \in \tE, i, \ell \in N_G,\\
    %& q_{i\ell} \in \{0,1\}, \;\;\forall  i, \ell \in N_G. \label{eq:a9}
%\end{align}
% ---
 \begin{align}
    \min &\sum_{i, \ell \in N_G:\atop \{u_i,u_\ell\} \in E} d_{i\ell} q_{i\ell} \label{eq:a1}\\
    \mbox{s. t.} \;
    % ---
    &\sum_{i \in \RR_j^G} t_{ji} = 1, &&\forall j \in N, \label{eq:a2}\\
    % ---
    &\sum_{\ell \in N_G:\atop \{u_i,u_\ell\}\in U} \!\!\!\! f^{jk}_{i\ell} -
     \sum_{\ell \in N_G:\atop \{u_\ell,u_i\}\in U} \!\!\!\! f^{jk}_{\ell i} = t_{j i},
     &&\forall j, k \in N: (v_j, v_k) \in A, i \in N_G, \label{eq:a3}\\ 
    % ---
     &\sum_{\ell \in N_G:\atop \{u_\ell,u_i\}\in U} \!\!\!\!f^{jk}_{\ell i} - \sum_{\ell \in N_G:\atop \{u_i,u_\ell\}\in U}\!\!\!\! f^{jk}_{i\ell} = t_{ki},
     &&\forall j, k \in N: (v_j, v_k) \in A, i \in N_G, \label{eq:a4}\\
     & \dsum_{i, \ell \in N_G:\atop \{u_i,u_\ell\} \in E} (f_{i\ell}^{jk} + f_{\ell i}^{jk}) \leq q_{i\ell},  && \forall \, j, k \in N: (v_j,v_k) \in A,\label{eq:a6} \\
    % % % -------------------------------------
      &q_{i\ell} \geq f^{jk}_{i\ell}, 
      &&\forall i, \ell \in N_G: \{u_i,u_\ell\} \in E, \notag\\
      &\; &&\;\; j, k \in N': (v_j,v_k) \in A, \label{eq:a7} \\
      & f_{i\ell}^{jk} \geq 0, &&\forall j, k \in N': (v_j, v_k) \in A, i, \ell \in N_G,\\
    % ---
    & q_{i\ell} \in \{0,1\}, &&\forall  i, \ell \in N_G. \label{eq:a9}
\end{align}
The objective function accounts for the total length of the tree. Constraints \eqref{eq:a2} enforce that all nodes in $V$ are assigned to exactly one node in $U$ (the embedding). Constraints \eqref{eq:a3} and \eqref{eq:a4} are the path-based constraints that ensure that a unit of flow is sent between the embeddings of nodes $j$ and $k$. Constraints \eqref{eq:a6} avoid the use of both direction of the edges in the solution and that different origin-destination pairs share the same edges in the solution. Constraints \eqref{eq:a7} correctly define the variable $q$.

\begin{figure}[h] 
\begin{center}
\includegraphics[trim={0.2cm 0.2cm 0.2cm 0.2cm},clip,width=1\textwidth]{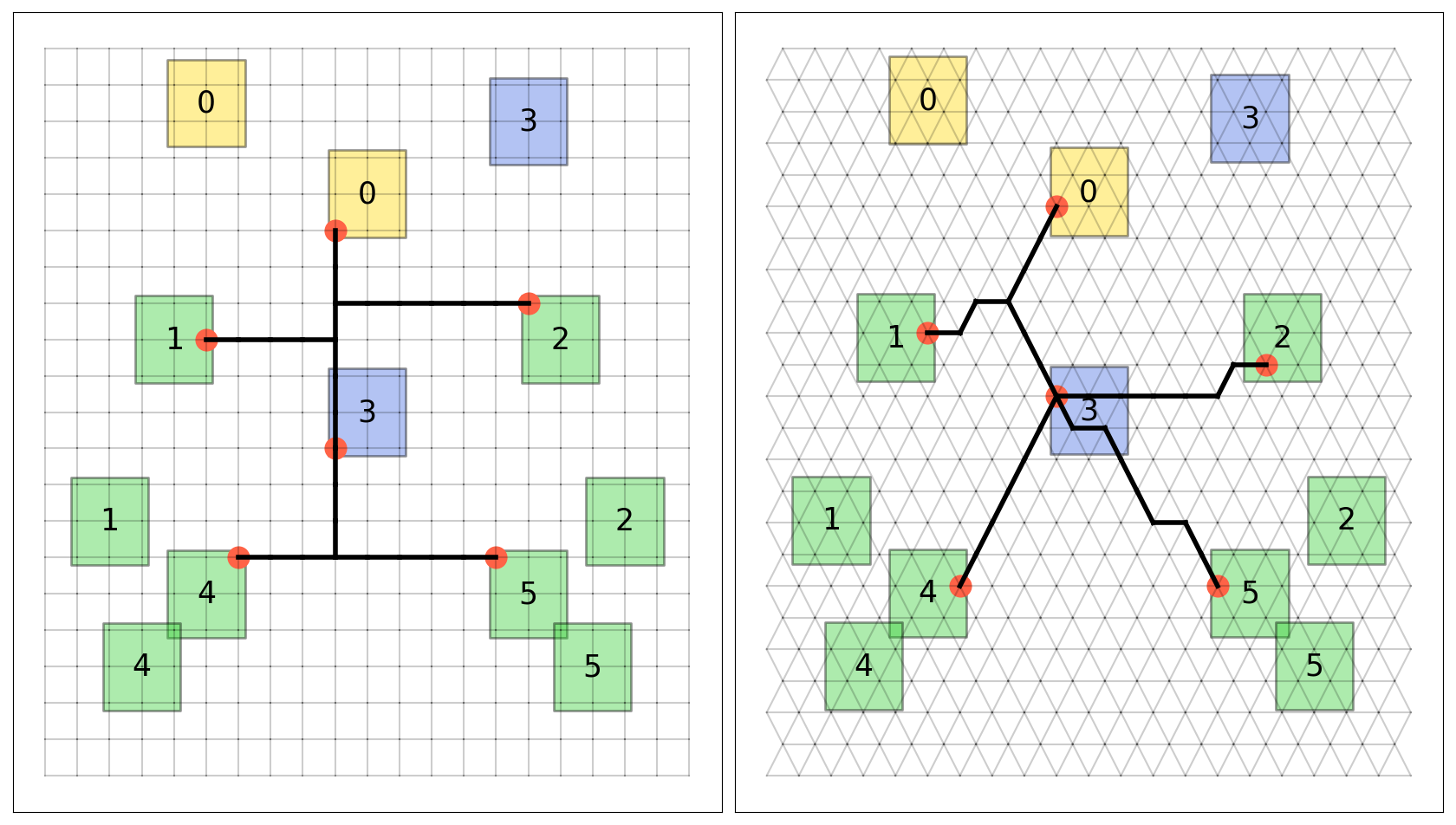}
\end{center}
\caption{Graphical representation of the solution for scenario $\ref{fig:totally_unsolved}$ using the formulation (\ref{eq:a1}-\ref{eq:a9}) for two different graphs $G$. \label{fig:kkk}}
\end{figure}

In Figure \ref{fig:kkk}, we observe two distinct discretizations of the example scenario \ref{fig:totally_unsolved}, along with the corresponding solutions obtained by solving the model (\ref{eq:a1}-\ref{eq:a9}). On the left, the solution space has been discretized using a square grid, while on the right, a triangular grid is utilized. Additionally, the positions of the final nodes are highlighted in orange, and the Steiner tree, representing the solution to the problem, is depicted in black. 

% --------------------------------------------------
% --------------------------------------------------
% --------------------------------------------------
% --------------------------------------------------
% --------------------------------------------------
% --------------------------------------------------
% --------------------------------------------------
% --------------------------------------------------
% --------------------------------------------------
% --------------------------------------------------
% --------------------------------------------------
% --------------------------------------------------
% --------------------------------------------------
% --------------------------------------------------
% --------------------------------------------------

\subsection{A variable fixing strategy for geometric graphs \label{subsec:improving}}

The graph $G$, representing the solution space in the discrete-domain case, can be dense, particularly when it is generated to discretize a continuous space. In this context, the graph $G$ is typically a geometric graph inducing a grid in $\R^d$. In Figure \ref{fig:kkk} we showed dense grids (square and triangular) used as graph $G$. In such graphs, the special hierarchical shape of the problem may help to fix some of the $f$ and $q$ variables to zero, reducing the problem size and enabling the solution of larger instances more quickly. This reduction is more efficiently applied by removing parts of the graph that the paths between nodes will never use.

Let $v\in V'$ and $\varrho^{-1}(v)$ its set of children. The paths in $G$ linking $v$ with its children start at some point in $\RR_v$ and end at points in $\bigcup_{w \in \varrho^{-1}(v)} \RR_w$. Let $\mathcal{C}_v$ be the convex (for $\ell_p$-norms with $1<p<\infty$) or rectangular (for $p=1, \infty$) hull of these two sets. Thus, the paths linking $v$ with its children in the optimal solution, since the lengths of the edges are induced by norm-based distances will never use edges outside $\mathcal{C}_v \cap E$. Consequently, the variables $f_{i\ell}^{jk}$ and $q_{i\ell}$ can be fixed to zero for any $\{i,\ell\} \in E$ that is not in $\mathcal{C}_v$ ($j$ is the index identified with $v$ and $k$ the indices of each of its children).

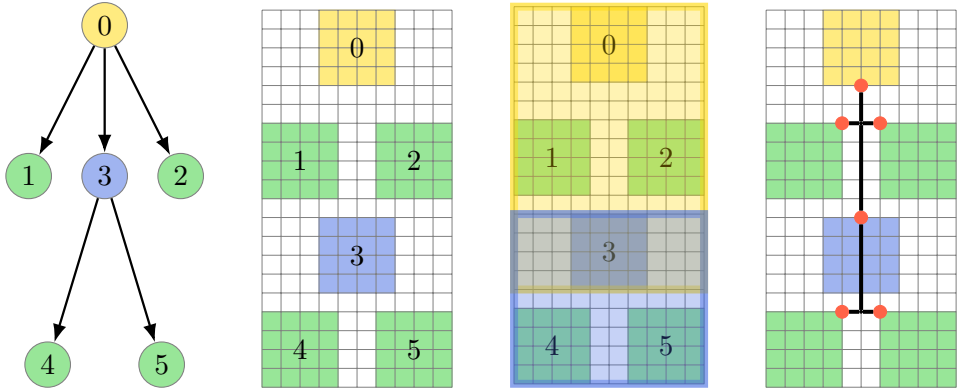
\begin{figure}[h]
\begin{tikzpicture}[scale=0.25]
\node (A0) at (5,20) [circle, draw, fill=gold, inner sep=1, opacity=.5, text opacity=1, minimum size=6mm] {\normalsize $0$};
\node (A1) at (1,12) [circle, draw, fill=limegreen, inner sep=1, opacity=.5, text opacity=1, minimum size=6mm] {\normalsize $1$};
\node (A2) at (9,12) [circle, draw, fill=limegreen, inner sep=1, opacity=.5, text opacity=1, minimum size=6mm] {\normalsize $2$};
\node (A3) at (5,12) [circle, draw, fill=royalblue, inner sep=1, opacity=.5, text opacity=1, minimum size=6mm] {\normalsize $3$};
\node (A4) at (2,2) [circle, draw, fill=limegreen, inner sep=1, opacity=.5, text opacity=1, minimum size=6mm] {\normalsize  $4$};
\node (A5) at (8,2) [circle, draw, fill=limegreen, inner sep=1, opacity=.5, text opacity=1, minimum size=6mm] {\normalsize $5$};
\draw[-Latex, line width=0.3mm] (A0)--(A1);
\draw[-Latex, line width=0.3mm] (A0)--(A2);
\draw[-Latex, line width=0.3mm] (A0)--(A3);
\draw[-Latex, line width=0.3mm] (A3)--(A4);
\draw[-Latex, line width=0.3mm] (A3)--(A5);
\end{tikzpicture}~\hspace*{0.5cm}
% -----------------------------
\begin{tikzpicture}[scale=0.25]
\draw[draw=gold, fill=gold, opacity=0.5] (3,16) rectangle ++(4,4);
\draw[draw=limegreen, fill=limegreen, opacity=0.5] (0,10) rectangle ++(4,4);
\draw[draw=limegreen, fill=limegreen, opacity=0.5] (6,10) rectangle ++(4,4);
\draw[draw=royalblue, fill=royalblue, opacity=0.5] (3,5) rectangle ++(4,4);
\draw[draw=limegreen, fill=limegreen, opacity=0.5] (0,0) rectangle ++(4,4);
\draw[draw=limegreen, fill=limegreen, opacity=0.5] (6,0) rectangle ++(4,4);
\draw [step=1.0, very thin, gray] (0,0) grid (10,20); 
\node (A0) at (5,18) {$0$};
\node (A1) at (2,12) {$1$};
\node (A2) at (8,12) {$2$};
\node (A3) at (5,7) {$3$};
\node (A4) at (2,2) {$4$};
\node (A5) at (8,2) {$5$};
\end{tikzpicture}~\hspace*{0.5cm}
% -----------------------------
\begin{tikzpicture}[scale=0.25]
\draw[draw=gold, fill=gold, opacity=0.5] (3,16) rectangle ++(4,4);
\draw[draw=limegreen, fill=limegreen, opacity=0.5] (0,10) rectangle ++(4,4);
\draw[draw=limegreen, fill=limegreen, opacity=0.5] (6,10) rectangle ++(4,4);
\draw[draw=royalblue, fill=royalblue, opacity=0.5] (3,5) rectangle ++(4,4);
\draw[draw=limegreen, fill=limegreen, opacity=0.5] (0,0) rectangle ++(4,4);
\draw[draw=limegreen, fill=limegreen, opacity=0.5] (6,0) rectangle ++(4,4);
\draw [step=1.0, very thin, gray] (0,0) grid (10,20); 
\node (A0) at (5,18) {$0$};
\node (A1) at (2,12) {$1$};
\node (A2) at (8,12) {$2$};
\node (A3) at (5,7) {$3$};
\node (A4) at (2,2) {$4$};
\node (A5) at (8,2) {$5$};
\draw[draw=gold, line width=3pt, fill=gold, opacity = 0.6, fill opacity=0.3] (0,5) rectangle ++(10,15);
\draw[draw=royalblue, line width=3pt, fill=royalblue, opacity = 0.6, fill opacity=0.3] (0,0) rectangle ++(10,9);
\end{tikzpicture}~\hspace*{0.5cm}
% -----------------------------
\begin{tikzpicture}[scale=0.25]
\tikzset{myarrow/.style={decoration={markings,mark=at position 1 with %
    {\arrow[scale=0.75,>=stealth]{>}}},postaction={decorate}}}
\draw[draw=gold, fill=gold, opacity=0.5] (3,16) rectangle ++(4,4);
\draw[draw=limegreen, fill=limegreen, opacity=0.5] (0,10) rectangle ++(4,4);
\draw[draw=limegreen, fill=limegreen, opacity=0.5] (6,10) rectangle ++(4,4);
\draw[draw=royalblue, fill=royalblue, opacity=0.5] (3,5) rectangle ++(4,4);
\draw[draw=limegreen, fill=limegreen, opacity=0.5] (0,0) rectangle ++(4,4);
\draw[draw=limegreen, fill=limegreen, opacity=0.5] (6,0) rectangle ++(4,4);
\draw [step=1.0, very thin, gray] (0,0) grid (10,20); 
\node (A0) at (5,18) {};
\node (A1) at (2,12) {};
\node (A2) at (8,12) {};
\node (A3) at (5,7) {};
\node (A4) at (2,2) {};
\node (A5) at (8,2) {};
\node(B0) at (5,16) [draw, fill, circle, inner sep=1.7, color=tomato] {};
%\node(S1) at (5,14) [draw, diamond, inner sep=1] {};
\node(S1) at (5,14) [draw, fill, circle, inner sep=0, opacity=1] {};
\node(B1) at (4,14) [draw, fill, circle, inner sep=1.7, color=tomato] {};
\node(B2) at (6,14) [draw, fill, circle, inner sep=1.7, color=tomato] {};
\node(B3) at (5,9) [draw, fill, circle, inner sep=1.7, color=tomato] {};
%\node(S2) at (5,4) [draw, diamond, inner sep=1] {};
\node(S2) at (5,4) [draw, fill, circle, inner sep=0, opacity=1] {};
\node(B4) at (4,4) [draw, fill, circle, inner sep=1.7, color=tomato] {};
\node(B5) at (6,4) [draw, fill, circle, inner sep=1.7, color=tomato] {};
\draw[line width=1.5pt] (B0)--(S1);
\draw[line width=1.5pt] (S1)--(B1);
\draw[line width=1.5pt] (S1)--(B2);
\draw[line width=1.5pt] (S1)--(B3);
\draw[line width=1.5pt] (B3)--(S2);
\draw[line width=1.5pt] (S2)--(B4);
\draw[line width=1.5pt] (S2)--(B5);
\end{tikzpicture}
\caption{Illustration of the variable fixing strategy for geometric graphs.\label{f:varfix}}
\end{figure}

In Figure \ref{f:varfix}, we illustrate this construction for a toy instance. In the two left plots, we show the topology of the tree $T$ and the distribution of the neighborhoods in the grid graph, respectively. In the third picture, we highlight with a yellow box the space in the graph where the paths from $0$ to $1$, $2$, and $3$ will be traced, and in blue the space where the paths from $3$ to $4$ and $5$ will be routed in an optimal solution. Thus, all the $f$ variables for nodes/edges outside the yellow rectangle can be fixed to $0$ for those links and assignments of the first depth of the tree, and analogously for the blue box. In the right plot, we show an optimal solution for this instance.

\section{Experiments \label{sec:experiments}}

In this section, we present the results of our computational experiments, conducted to validate the performance of our approaches. The specific instances that we generate as well as the detailed results are available at 
\href{https://github.com/anticiclon/Fixed-Topology-Minimum-Length-Trees-with-Neighborhoods}{\small\texttt{github.com/anticiclon/Fixed-Topology-Minimum-Length-Trees-with-Neighborhoods}}.

\subsection{Instances}

We generated a series of random instances with different characteristics. Each instance consists of a random tree with $n$ nodes, with $n \in \{20, 50, 100, 200\}$ using the \texttt{random\_tree} function in the python library \texttt{networkx}. The generated trees are relabeled so that node $0$ is the root node, and it is directed from root to leaves.

We embed the tree in the square $[0,2000] \times [0,2000]$. For each tree, we generate $b$ planar unions of square neighborhoods for $b \in \{1,2,3,4,5\}$ with different sizes $s \in \{20, 50, 100, 200\}$. Combining the different values for $n$, $b$, and $s$, the total number of instances generated was $400$.

\begin{figure}[h] 
\begin{center}
\includegraphics[trim={0cm 0cm 0cm 0cm}, clip, width=1\textwidth]{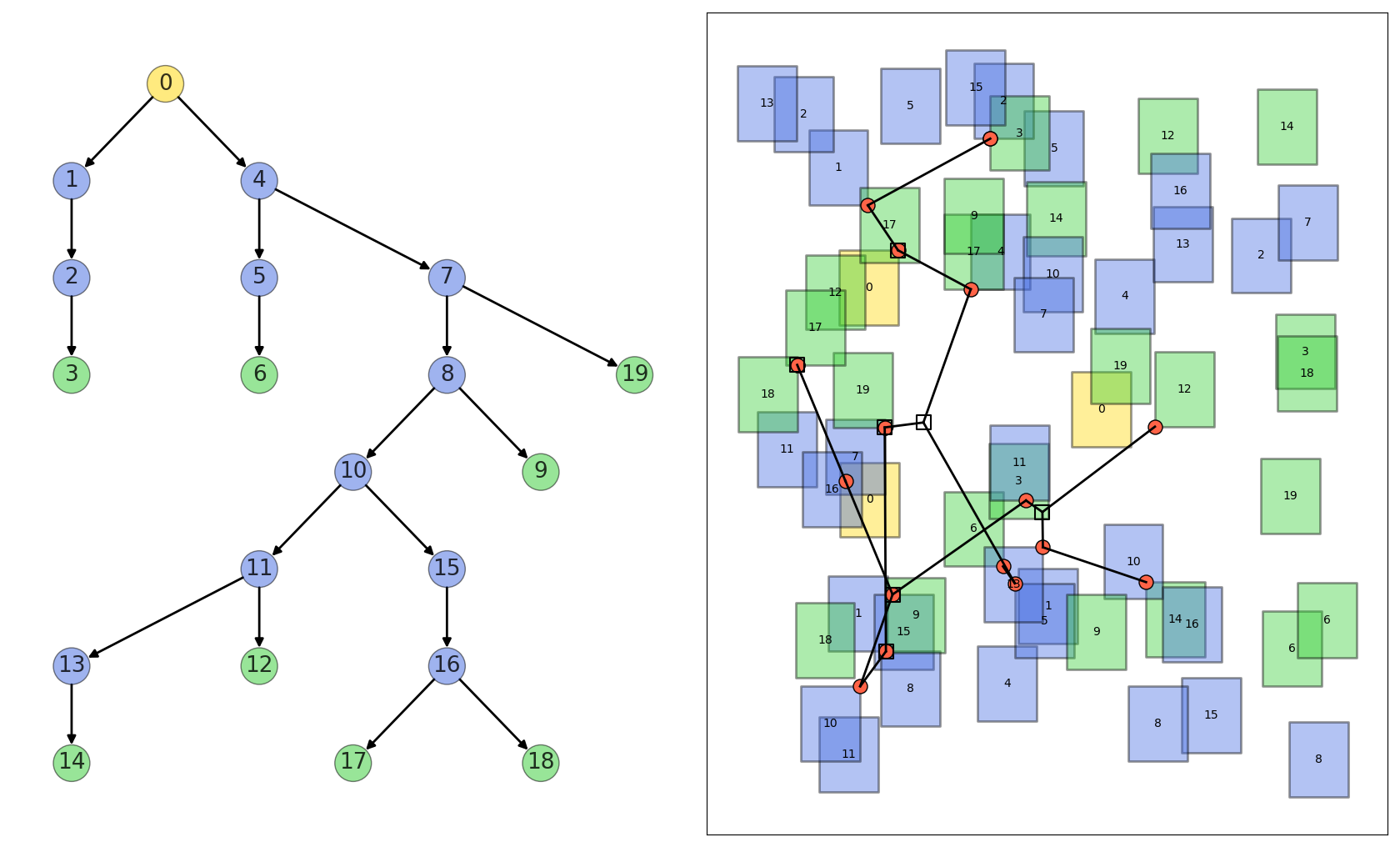}
\end{center}
\caption{An example of a random instance with 20 nodes and 3 spaces per
node solved by the continuous model with norm $\ell_{2}$. \label{fig:yyy}}
\end{figure}

In Figure \ref{fig:yyy}, we present an example of a random instance with $n=20$ nodes, $b=3$ connected components for the neighborhoods, and square side length $s=200$, as well as the solution for the Euclidean norm.

All these scenarios were solved using three models:

\begin{description}
    \item[Cont $\ell_1$:]  Continuous-domain FTTPN with $\ell_1$-norm based distances.
        \item[Cont $\ell_2$:]  Continuous-domain FTTPN with $\ell_2$-norm based distances.
    \item[Disc:] We construct ad-hoc graphs $G$ based on each instance. For parent-children combinations, we generate a rectangular grid by projecting the vertices of the neighborhoods of parents and children and tracing orthogonal lines. 
\end{description}

All instances were solved using the Gurobi 10.0.3 optimizer on an Ubuntu 22.04.3 environment, with an AMD EPYC 7042p 24-Core Processor and 64 GB RAM. Initially, all parameters of the Gurobi solver were set to their default values, and a CPU time limit of 7200 seconds was set.

A first tabular summary of the results is reported  in Table \ref{table:1}. In that table, for each of the approaches, for each node size ($n$) and each number of components for the neighborhoods ($b$), we report the average required CPU time to solve the instances (only those up to the $20$ per row that were optimally solved), the MIPGAP (of those not optimally solved), and the number of non-optimaly solved instances (out of $20$). For the {\bf Disc} approach, all the instances were optimally solved, so we do not report the MIPGAP (all are $0\%$) nor the number of unsolved instances. The flag TL indicates that none of the instances summarized in that row could be optimally solved within the time limit (and then, the MIPgap is the average gap of all the instances in the row).

\begin{table}[h!]
    \centering
    \begin{adjustbox}{max width=\textwidth}
    \begin{tabular}{|c|c|c|c|c|c|c|c|c|}
        \cline{3-9}
        \multicolumn{2}{c}{} & \multicolumn{3}{|c}{{\bf Cont} $\ell_1$} & \multicolumn{3}{|c}{{\bf Cont} $\ell_2$} & \multicolumn{1}{|c|}{{\bf Disc}}\\\hline
$n$ & $b$ & CPUTime & MIPGAP & UnSolved & CPUTime & MIPGAP & UnSolved & CPUTime \\ \hline
        \multirow{5}{*}{20} & 1 & 7.29 & 0\% & 0 & 30.04 & 0\% & 0 & 0.06 \\ 
                             & 2 & 0.58 & 0\% & 0 & 3.91 & 0\% & 0 & 0.34 \\ 
                             & 3 & 4.18 & 0\% & 0 & 62.05 & 0\% & 0 & 1.61 \\ 
                             & 4 & 7.72 & 0\% & 0 & 589.73 & 0\% & 0 & 6.43 \\ 
                             & 5 & 17.16 & 0\% & 0 & 361.46 & 13.9\% & 2 & 15.49 \\ \hline
        \multirow{5}{*}{50} & 1 & 81.20 & 4.99\% & 1 & 370.49 & 11.7\% & 3 & 0.19 \\ 
                             & 2 & 227.84 & 22.68\% & 6 & 1301.79 & 7.4\% & 5 & 2.14 \\ 
                             & 3 & 354.74 & 13.52\% & 11 & 1795.31 & 8.7\% & 10 & 6.33 \\ 
                             & 4 & 1157.23 & 0\% & 0 & 3330.50 & 18.7\% & 17 & 21.91 \\ 
                             & 5 & 3613.66 & 0\% & 0 & TL & 22.7\% & 20 & 20.76 \\ \hline
        \multirow{5}{*}{100} & 1 & 214.99 & 0\% & 0 & 1822.62 & 5.4\% & 16 & 0.48 \\ 
                              & 2 & 2723.87 & 7.07\% & 8 & 16.75 & 9.6\% & 17 & 6.61 \\ 
                              & 3 & 1527.23 & 8.92\% & 17 & 20.09 & 15.5\% & 19 & 17.23 \\ 
                              & 4 & TL & 22.30\% & 20 & TL & 30.7\% & 20 & 28.93 \\ 
                              & 5 & TL & 31.91\% & 20 & TL & 39.5\% & 20 & 105.03 \\ \hline
        \multirow{5}{*}{200} & 1 & 123.87 & 2.76\% & 1 & 16.71 & 7.7\% & 20 & 1.25 \\ 
                              & 2 & 684.94 & 8.44\% & 17 & 20.11 & 11.9\% & 17 & 14.96 \\ 
                              & 3 & TL & 16.52\% & 20 & TL & 23.2\% & 18 & 57.38 \\ 
                              & 4 & TL & 27.37\% & 20 & TL & 35.4\% & 20 & 65.95 \\ 
                              & 5 & TL & 37.74\% & 20 & TL & 47.3\% & 20 & 237.23 \\ \hline
    \end{tabular}
    \end{adjustbox}
    \caption{Average results of our computational experiments. \label{table:1}}
    \label{tab:performance}
\end{table}
As expected, one can observe that as the number of nodes increases, the problem becomes more challenging. This also happens with the number of components for the neighborhoods ($b$).

The discrete approach was able to find solutions much faster than the continuous approaches.

\begin{figure}[h] 
\begin{center}
\includegraphics[width=0.8\textwidth]{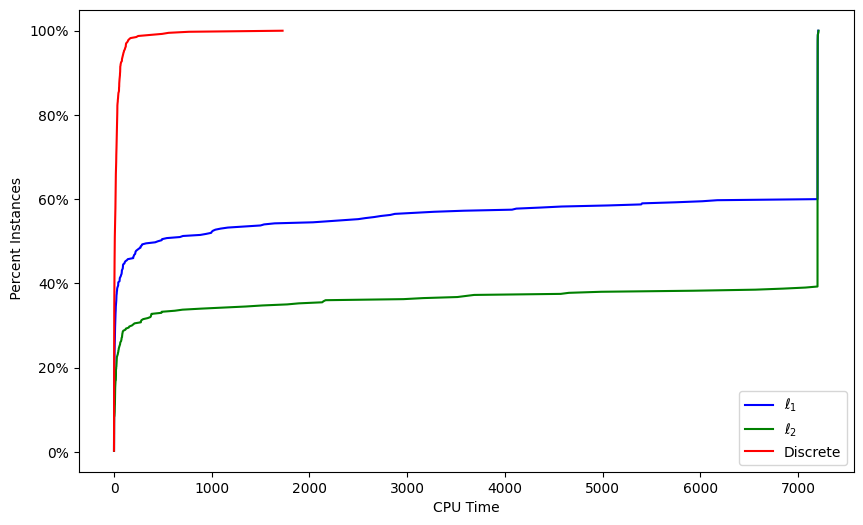}
\end{center}
\caption{Graphical representation of the time-based performance profiles for the three methods. \label{fig:pp}}
\end{figure}

In Figure \ref{fig:pp} we plot the performance profiles for the three methods on all the instances. On the $x$-axis we represent the CPU time. On the $y$-axis, the percentage of instances that were optimally solved within a time smaller than the one on the $x$-axis is shown. As already observed in the table, the discrete approach solved the instances in less than $2000$ seconds, with most of them ($80\%$) solved in less than $100$ seconds. Among the continuous approaches, the $\ell_2$ approach is more challenging, with only $40\%$ of the instances optimally solved within the time limit. The $\ell_1$ norm approach (whose formulation is a Mixed Integer Linear problem) solved $60\%$ of the instances within the time limit.

\begin{figure}[h] 
\begin{center}
\includegraphics[width=0.8\textwidth]{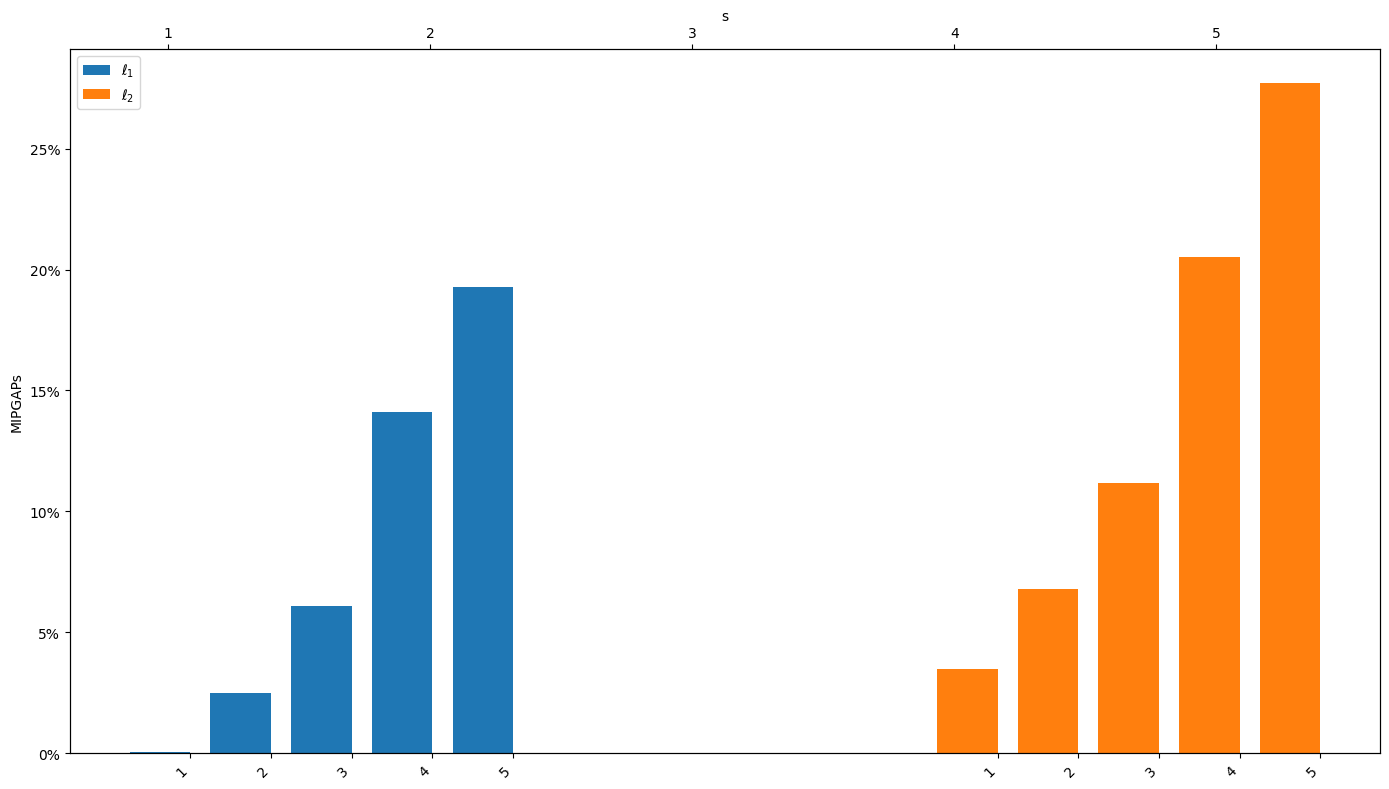}
\end{center}
\caption{Graphical representation of MIPGAPs for the continuous approaches. \label{fig:gaps}}
\end{figure}

We show in Figure \ref{fig:gaps} the MIPGAPs aggregated by the number of components in the neighborhoods. One can observe that the MIPGAP clearly increases with the number of components. Note that the components for each node are modeled as disjunctive constraints, introducing binary variables for each component. Therefore, as the number of components increases, the problem becomes more difficult.

\begin{figure}[h] 
\begin{center}
\includegraphics[width=0.8\textwidth]{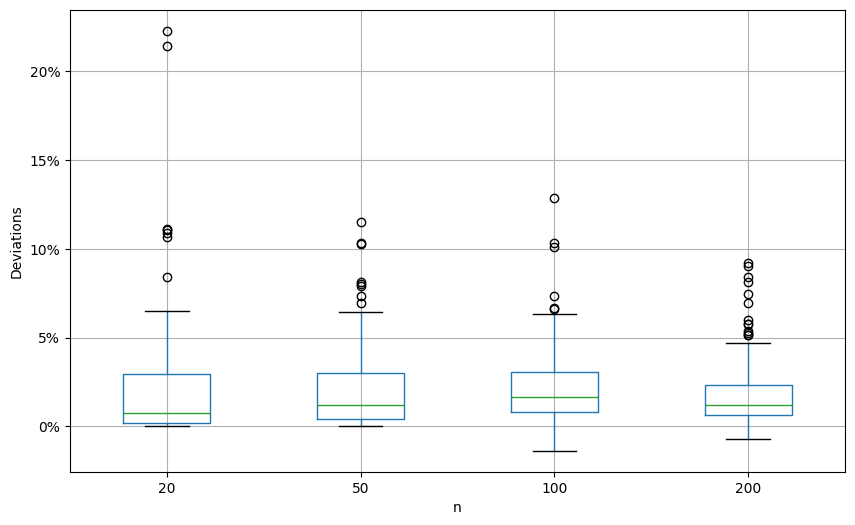}
\end{center}
\caption{Graphical comparison of the objective functions of both models: the Discrete model with grid graph and the continuous model with the $\ell_{1}$ norm.\label{fig:devs}}
\end{figure}

Finally, since we construct the graph for the Discrete approach using a grid, we analyze how close the solution of the Discrete approach is to the solution of the continuous $\ell_1$. In Figure \ref{fig:devs} we plot a boxplot showing the deviation of the objective value of the Discrete approach from the objective value of the continuous $\ell_1$ approach, for each size($s$) considered in our instances. The deviation is measured as:
$$
{\rm Dev} = \frac{\text{objective value {\bf Disc}} -\text{objective value {\bf Cont} $\ell_1$}}{\text{objective value {\bf Disc}}} \times 100
$$
Note that the deviations are very small, most of them close to zero. Thus, the discrete approach appears to be a promising method to approximate the solution of the Continuous $\ell_1$ problem, solving the instances in much less CPU time.

To summarize, the continuous model with $\ell_{2}$ distance produces the shortest trees but is also the slowest, especially for larger instances. The continuous model with $\ell_{1}$ distance performs well for small to mid-sized instances but also struggles with larger instances. The discrete model is the most efficient in terms of computational time. For large instances, all models face challenges with computational time, indicating a need for further optimization or alternative approaches to handle large-scale problems efficiently.

% --------------------------------------------------
% --------------------------------------------------
% --------------------------------------------------
% --------------------------------------------------
% --------------------------------------------------
% --------------------------------------------------
% --------------------------------------------------
% --------------------------------------------------
% --------------------------------------------------
% --------------------------------------------------
% --------------------------------------------------
% --------------------------------------------------
% --------------------------------------------------
% --------------------------------------------------
% --------------------------------------------------

\section{Conclusions \label{sec:conclusions}}

In this paper, we introduce a new problem, the FTTNP, motivated by its various industrial applications. We provide a general framework for the problem and study two versions of practical interest, the continuous and the discrete-domain frameworks. We propose mathematical optimization approaches to tackle both versions of the problem. Our approach combines different elements: from the transformation of the given instance to the analysis of the problem as a minimum-length Steiner tree problem. For the continuous case, we develop a formulation that can be rewritten as a Mixed Integer Second Order Cone problem, as well as some simplifications for specific instances of interest. For the discrete FTSTNP, we derive a mixed-integer linear programming approach. 

Several extensions are possible for futher research on this problem. On the one hand, in this type of pipeline/cables routing problems that arise in industrial applications it is usually required to explicitly account for the use of \textit{elbows} in the routes (in some applications they are considerably costly). Thus, the graph where the cables are traced must consider changes of directions of the routes. In \cite{blanco2021network}, the authors propose a transformation of the graph, together with a modification of the objective function, that allows this specification. Its adaptation to our problem, and its implications in terms of complexity of the problems have to be studied. On the other hand, the development of heuristic approaches for the problem is advisable since in some applications the size of the discrete optimization problems that we propose can be large, being the available optimization solvers unable to solve them. The design of decomposition strategies can be good option for this problem

% --------------------------------------------------
% --------------------------------------------------
% --------------------------------------------------
% --------------------------------------------------
% --------------------------------------------------
% --------------------------------------------------
% --------------------------------------------------
% --------------------------------------------------
% --------------------------------------------------
% --------------------------------------------------
% --------------------------------------------------
% --------------------------------------------------
% --------------------------------------------------
% --------------------------------------------------
% --------------------------------------------------
\section*{Acknowledgements}

This research has been partially supported by grant PID2020-114594GB-C21 funded by MICIU/AEI/10.13039/501100011033;  grant RED2022-134149-T funded by MICIU/AEI
/10.13039/501100011033 (Thematic Network on Location Science and Related Problems); grant C-EXP-139-UGR23 funded by the Consejería de Universidad,
Investigación e Innovación and by the ERDF Andalusia Program 2021-2027, grant AT 21\_00032, and the IMAG-María de Maeztu grant CEX2020-001105-M /AEI /10.13039/501100011033.

\end{document}